\documentclass[a4paper,10 pt,titlepage]{report}
\usepackage[english]{babel}
\usepackage{amsthm}
\usepackage{graphicx}
\usepackage{geometry}
\usepackage{authblk}

\usepackage{amsmath,amsfonts}
\theoremstyle{definition}

\theoremstyle{plain}

\title{Computational issues and numerical experiments for Linear Multistep Method Particle Filtering}

\author[1]{Daniela Calvetti A\thanks{dxc57@case.edu}}
\author[2]{Salvatore Cuomo B\thanks{salvatore.cuomo@unina.it}}
\author[2]{Monica Pragliola\thanks{m.pragliola@studenti.unina.it}}
\author[1]{Erkki Somersalo\thanks{ejs49@case.edu}}
\author[2]{Gerardo Toraldo\thanks{toraldo@unina.it}}
\affil[1]{Department of Mathematics, Case Western Reserve University, 10900 Euclid Ave., Cleveland, OH 44106, USA}
\affil[2]{Department of Mathematics and Applications, University of Naples Federico II, Strada Vicinale Cupa Cintia 21, 80126, Naples, Italy.}
\date{}
\begin{document}
	\maketitle
\section{Introduction}
The Linear Multistep Method Particle Filter (LMM PF) is a method for predicting the evolution in time of a evolutionary system governed by a system of differential equations. In order to appreciate the contribution of the LMM PF, let us suppose to deal with an ODEs system modelling an inverse problem and depending on unknown or poorly known parameters. The estimate of states and parameters for such a system from noisy measurements of a function of some of the states at discrete times is a central problem in several applications. In general, inverse problems are ill posed, that means, for instance, that the solution does not exist.\\
A possible way to overcome the mentioned difficulties is to resort to a statistical approach, outlining a framework in which the unknown states and parameters are modelled as random variables and the uncertainties on them are represented by probability density functions.
A great advantage in adopting such a point of view is that we do not have to worry about the existence of the solution, or of a unique solution, since the solution is not a single value, but a probability density function (posterior density). Nevertheless the original problem is not formulated in statistical terms. To find the value of the generic unknown, we can, for example, maximize the density function of the corresponding random variable. However, this optimization problem can not be solved directly. We need to appeal to a pre-processing phase that makes the issue practicable from a computational point of view, resorting to \textit{sampling} techniques, such as LMM PF. The aim of sampling techniques is to draw information from probability density functions, whose analytical form is approximated by sample points. The computational efficiency of the sampling technique is crucial for the success of the method, since, the better the sampling, the more accurate and reliable is the final solution. \\ 
Let us present a more detailed analysis of the issue and consider a system of differential equations depending on a vector of unknown parameters $\theta$:
\begin{equation}
\frac{du}{dt}=f(t,u,\theta),\quad u(0)=u_{0},\quad t\in{[0,T]}.
\end{equation}
where $u = u(t)\in R^{d}$ is a vector containing the states of the system, $f:\mathbb{R} \times \mathbb{R}^{d}\times \mathbb{R}^{k}\longrightarrow R^{d}$ is the known model function, and $\theta \in \mathbb{R}^{k}$  is the vector of model parameters. \\
Observe that (1) is not required to model an inverse problem, since the method proposed in the following can be successfully applied to every kind of ODEs system. 
Setting a discretization step for the time interval $[0,T]$, suppose that the measured observations are given by
\begin{equation}
b_{j}=g(u(t_{j}),\theta)+e_{j} \in \mathbb{R}^{m},
\end{equation}
where $g:\mathbb{R}^{d} \times \mathbb{R}^{k} \longrightarrow \mathbb{R}^{m}$ is a known function and $e_{j}$ denotes the noise in the measurement process, which is additive for sake of simplicity.\\
Formally, we are looking for an estimate of $u(t)$ at given times and $\theta$ from the measurements $b_{j}$. \\
Let us denote by $D_{j}$ the set of data accumulated up to time $t=t_{j}$,
\begin{equation*}
D_{j}=\{b_{1},...,b_{j}\}.
\end{equation*}
As mentioned before, adopting a statistical approach, the final solution is a probability density function. In particular, here we are interested in updating the posterior density from one time instant to the next:
\begin{equation}
\pi(u_{j},\theta|D_{j})\longrightarrow\pi(u_{j+1},\theta|D_{j+1}),
\end{equation}
where $u_{j}$ denotes the discrete approximation of the state vector $u(t_{j})$.\\
As we are working in a statistical framework, the analytical model (1) must be converted into a statistical one, the so called \textit{evolution-observation model}. Let us show a possible way to carry out the mentioned transformation.
\\Consider (1) over the time interval $[t_{j},t_{j+1}]$
\begin{equation*}
\frac{du}{dt}=f(t,u,\theta), \quad t_{j}<t<t_{j+1},
\end{equation*}
and let $\psi^{exact}$ be the formal exact propagation operator:
\begin{equation*}
u(t_{j+1})=\psi^{exact}(t_{j+1},u(t_{j},\theta)).
\end{equation*}
We need to replace $\psi^{exact}$ by a numerical scheme. In particular, we choose an $r$-step solver $\psi$ such that
\begin{equation}
u_{j+1}=\psi(u_{j},u_{j-1},...,u_{j-r+1},\theta,h),
\end{equation}
where $h$ is the constant time step.
Substituting in (4) the exact solution for the numerical one, the equality is retained if we take into account the approximation error, or local truncation error:
\begin{equation*}
u(t_{j+1})=\psi(u(t_{j}),u(t_{j-1}),...,u(t_{j-r+1}),\theta,h)+w_{j+1}.
\end{equation*}
If we look at $u(t_{j})$ as the realization of the random variable $U_{j}$, the previous formula defines an $r$-Markov model:
\begin{equation}
U_{j+1}=\psi(U_{j},U_{j-1},...,U_{j-r+1},\theta,h)+W_{j+1},
\end{equation}
for the stochastic process $\{U_{j}\}_{j=0}^{T}$, where $\{W_{j}\}_{j=1}^{T}$ is the associated innovation process.\\
It is possible to turn the $r$-Markov model (5) into a $1$-Markov model by a change of variables:
\[X_{j}=
\begin{bmatrix}
U_{j} \\
.\\
.\\
.\\
U_{j-r+1}
\end{bmatrix}
\quad
V_{j+1}=
\begin{bmatrix}
W_{j+1} \\
0\\
.\\
.\\
.\\
0
\end{bmatrix}
\quad
\Psi(X_{j},\theta,h)=
\begin{bmatrix}
\psi(X_{j},\theta,h) \\
U_{j}\\
.\\
.\\
.\\
U_{j-r+2}
\end{bmatrix}
\]
\\
In the end we obtain the $1$-Markov model:
\begin{equation}
X_{j+1}=\Psi(X_{j},\theta,h)+V_{j+1}.
\end{equation}
Assume that $\{Y_{j}\}_{j=1}^{T}$ is the stochastic process modelling the observations of $U_{j}$. Then we have
\begin{equation}
Y_{j}=G(X_{j})+E_{j},
\end{equation}
where $\{E_{j}\}_{j=1}^{T}$ represents the measurement noise.\\
Equations (6)-(7) constitute the discrete-time evolution-observation model obtained from a discrete-time propagation system. 
It is worth observing that the innovation term $V_{j+1}$ in (6) represents mainly the numerical approximation error due to the propagation scheme $\psi$, since most of its components are null.

\section{LMM PF} 
Let us consider the discrete-time evolution-observation model (6)-(7):
\begin{equation}
\begin{split}
X_{j+1}&=\Psi(X_{j},\theta,h)+V_{j+1}, \qquad V_{j+1}\sim \mathcal{N}(0,\Gamma_{j+1}(X_{j},\theta))\\
Y_{j}&=G(X_{j})+E_{j},\qquad E_{j}\sim \mathcal{N}(0,\Sigma_{j}).
\end{split}
\end{equation}
As mentioned in the previous section, $V_{j+1}$ is connected with the approximation error, in fact its covariance matrix $\Gamma_{j+1}$ is computed by resorting to error estimate strategies, such as the Higher Order Method Error Control strategy (HOMEC). Adopting HOMEC, we need to consider an LMM method of order $p$ and an LMM method of order $\hat{p}\geq p+1$ from the same family. Denoting with $u$ the solution computed by the LMM method of order $p$, and with $\hat{u}$ the solution computed by the LMM method of order $\hat{p}$, we obtain the following expression for the innovation covariance
$$\Gamma_{j+1}=diag(\gamma),\quad \gamma_{i}=\tau^{2}(u_{j+1}-\hat{u}_{j+1})_{i}^{2},\quad i=1,...,d,$$
where $\tau>1$ is introduced to compensate for the omission of the higher order terms.\\
Our purpose is to compute a sequential update of the posterior density such as (3):
\begin{equation}
\pi(x_{j},\theta|D_{j})\longrightarrow \pi(x_{j+1},\theta|D_{j+1}).
\end{equation}
As in the classical particle filter [ref], we deal with samples approximating densities. Hence, the update we are actually interested in is:
$$S_{j}\longrightarrow S_{j+1}, \quad S_{j}=\big\{(x_{j}^{n},\theta_{j}^{n},w_{j}^{n})\big\}_{n=1}^{N},$$
where the pairs $(x_{j}^{n},\theta_{j}^{n})$ for $n=1,...,N$, have been drawn from the posterior density $\pi(x_{j},\theta|D_{j})$ with relative probabilities $w_{j}^{n}$.
\\For simplicity, we first assume that the parameter $\theta$ is known and can therefore be dropped from the notation in (4.2) and the sample $S_{j}$ has the form:
$$S_{j}=\big\{(x_{j}^{n},w_{j}^{n})\big\}_{n=1}^{N}.$$
The model (6)-(7) satisfies the following Markov properties:
\begin{itemize}
	\item the state $X_{j+1}$ depends on the past data $D_{j}$ only through the previous state, i.e.   $$\pi(x_{j+1}|x_{j},D_{j})=\pi(x_{j+1}|x_{j});$$
	\item the observation $Y_{j+1}$ depends on the past only through the current state $X_{j+1}$, i.e. 
	$$\pi(y_{j+1}|x_{j+1},D_{j})=\pi(y_{j+1}|x_{j+1}).$$
\end{itemize} 
The  previous properties suggest us to rewrite the chain of update (9) adding an intermediate step:
\begin{equation*}
\pi(x_{j}|D_{j})\longrightarrow\pi(x_{j+1}|D_{j})\longrightarrow\pi(x_{j+1}|D_{j+1}),
\end{equation*}
Hence, two updating formulas must be derived. The first one is based on the Chapman-Kolmogorov formula and controls the evolution update:
\begin{align*}
\pi(x_{j+1}|D_{j})&=\int \pi(x_{j+1}|x_{j},D_{j})\pi(x_{j}|D_{j})dx_{j}\\
&=\int \pi(x_{j+1}|x_{j})\pi(x_{j}|D_{j})dx_{j}.
\end{align*}
The second updating formula is obtained by applying the Bayes' formula and the Monte Carlo approximation and controls the observation update:
\begin{equation}
\begin{split}
\pi(x_{j+1}|D_{j+1}) & \propto \pi(y_{j+1}|x_{j+1}) \int \pi(x_{j+1}|x_{j})\pi(x_{j}|D_{j})dx_{j}\\
& \propto \pi(y_{j+1}|x_{j+1}) \sum_{n=1}^{N}w_{j}^{n} \pi(x_{j+1}|x_{j}^{n}).
\end{split}
\end{equation}
Basing on the layered sampling procedure, outlined in [ref], a new proposal $\hat{x}_{j+1}^{n}$ is drawn from the density $\pi(x_{j+1}|x_{j}^{n})$, for $n=1,...,N$. After the evaluation of the likelihoods $\pi(y_{j+1}|\hat{x}_{j+1}^{n})$ for each proposal particle, the algorithm resamples according to the relative likelihood. A possible consequence of this approach is the thinning of the sample, which is due to the discarding of the particles with the lowest likelihoods, and can be avoided by resorting to  the auxiliary particle strategy [ref]. So, we choose as auxiliary particle $\mu_{j+1}^{n}=\overline{x}_{j+1}^{n}$ for each $n$, where $\overline{x}_{j+1}^{n}=\Psi(x_{j}^{n},h)$ is the expectation of $X_{j+1}$ conditioned that $X_{j}=x_{j}^{n}$. In particular, $\overline{x}_{j+1}^{n}$ is a predictor of the value of $X_{j+1}$ given the initial value $x_{j}^{n}$.\\
We can now rewrite the updating formula (10)
$$\pi(x_{j+1}|D_{j+1}) \propto \sum_{n=1}^{N} w_{j}^{n}\pi(y_{j+1}|\overline{x}_{j+1}^{n})\frac{\pi(y_{j+1}|x_{j+1})}{\pi(y_{j+1}|\overline{x}_{j+1}^{n})}\pi(x_{j+1}|x_{j}^{n}),$$
The above equation can be seen as a mixture model, where $g_{j+1}^{n}=w_{j}^{n}\pi(y_{j+1}|\overline{x}_{j+1}^{n})$
is said to be the \textit{fitness} of the $n$-th \textit{predictor}.
In the following we are giving the algorithm of PF LMM for state estimation [ref]:\\

	\vspace{2mm}
	\textbf{Algorithm 3: LMM PF for state estimation}\\
	
	\textbf{Input}: $\pi(x_{0})=\pi(x_{0}|D_{0})$ prior distribution.
	\begin{enumerate}
		\item [(i)]\textbf{Initialize}: Draw the particle sample from $\pi(x_{0})$,\\ $S_{0}=\{(x_{0}^{1},w_{0}^{1}),...,(x_{0}^{N},w_{0}^{N})\}.$
		Set $j=0$;
		\item [(ii)]\textbf{Propagation}: Compute the predictor using LMM, \\
		$\overline{x}_{j+1}^{n}=\Psi(x_{j}^{n},h), \quad 1\leq n\leq N;$
		\item [(iii)]\textbf{Survival of the fittest}: For each $n$
		\begin{itemize}
			\item [-]Compute the normalized fitness weights:
			\begin{equation*}
			g_{j+1}^{n}=w_{j}^{n}\pi(y_{j+1}|\overline{x}_{j+1}^{n}), \quad g_{j+1}^{n}\longleftarrow \frac{g_{j+1}^{n}}{\sum_{n}g_{j+1}^{n}};
			\end{equation*}
			\item[-] Draw indices with replacement $l_{n}\in \{1,...,N\}$ using probabilities $P(l_{n}=k)=g_{j+1}^{k}$;
			\item[-] Reshuffle:
			$x_{j}^{n}\longleftarrow x_{j}^{l_{n}}, \quad \overline{x}_{j}^{n}\longleftarrow \overline{x}_{j}^{l_{n}}, \quad 1 \leq n \leq N.$
		\end{itemize}
		\item [(iv)]\textbf{Innovation}: For each $n$
		\begin{itemize}
			\item [-] Using LMM error control, estimate $\Gamma_{j+1}^{n}$;
			\item [-]Draw $v_{j+1}^{n} \sim \mathcal{N}(0,\Gamma_{j+1}^{n})$;
			\item [-]Proliferate: 
			$x_{j+1}^{n}=\overline{x}_{j+1}^{n}+v_{j+1}^{n}.$
		\end{itemize}
		\item[(v)] \textbf{Weight updating}: For each $n$ compute 
		\begin{equation*}
		w_{j+1}^{n}=\frac{\pi(y_{j+1}|x_{j+1}^{n})}{\pi(y_{j+1}|\overline{x}_{j+1}^{n})}, \quad w_{j+1}^{n}\longleftarrow \frac{w_{j+1}^{n}}{\sum_{n}w_{j+1}^{n}};
		\end{equation*}
		\item [(vi)] if $j<T$, set $j=j+1$ and repeat from Step 2; otherwise, stop.
	\end{enumerate}
	\textbf{Output}: $S_{k+1}$, $k=0,...,T-1$.	
	\vspace{2mm}

\vspace{10mm}
The LMM PF can also face up with parameter estimation problem, as clearly shown in [ref].
Let us have a quick review of the types of errors we are dealing with. First of all, we must consider the error introduced in order to generate the initial particle ensemble $S_{0}$, i.e. the \textbf{initial variance} $V_{0}$ of the prior distribution $\pi_{0}$. Furthermore, at each time instant $t_{j}$, the \textbf{variance} $V_{j}$ of the sample $S_{j}$ and the \textbf{absolute error} $E_{j}$ can be computed. In particular, denoted with $y^{exact}(t_{j})$ the analytical solution of (1) in $t_{j}$ and with $y^{mean}(t_{j})$ the mean of the particle ensemble $S_{j}$, we have
\begin{equation*}
E_{j}=\lvert y^{mean}(t_{j})-y^{exact}(t_{j})\rvert.
\end{equation*}
The absolute error takes into account the contribution of the \textit{global truncation error} and of the \textit{round off error} at each time instant $t_{j}$.
Our aim is to study the behaviour of the variances $V_{j}$ and of the errors $E_{j}$ for different values of $V_{0}$, when different numerical integration methods are considered.

\section{Computational experiments}
The LMM PF is applied to the dynamics of the skeletal muscle metabolism, in order to approximate concentrations of some metabolites in the blood and in the tissue during an episode of ischemia; $30$ out of $39$ are concentrations of metabolites in the tissue. The data consist of noisy observations of eight metabolites in the blood, which are glucose, lactate, alanine, triglyceride, glycerine, free fatty acid, carbon dioxide, and oxygen. The measurements are collected at $11$ time instants.\\
It is worth doing some considerations about the nuemrical results. At first, we choose the method BDF1 of order $p=1$ to propagate the particle ensemble, and the method BDF2 to get an estimate of the error.

\begin{figure}
	\begin{center}
		\includegraphics[width=70mm]{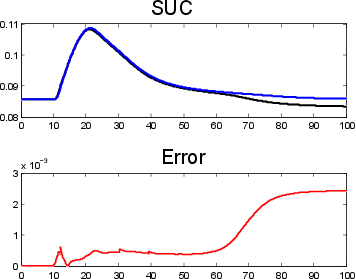}\\
		\includegraphics[width=70mm]{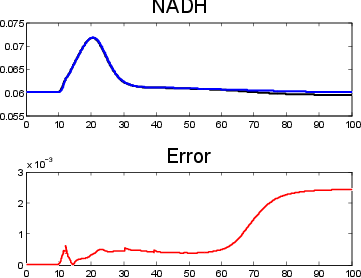}\\
		\includegraphics[width=70mm]{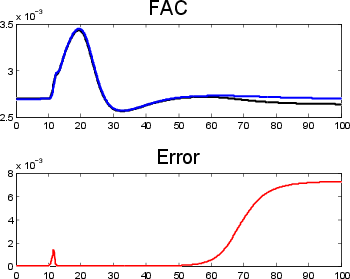}
	\end{center}
	\caption{Computed concentrations of succinate, fatty acyl-CoA and NAD reduced compared with their respective errors (BDF1-BDF2).}
\end{figure}

In Figure 1, the black solid line represents the approximated solution, while the blue solid line represents the real behavior of the metabolite taken into account. The red solid line describes the evolution of the absolute error, that is the difference between the approximated and the real concentration for each metabolite. We can observe that the error presents an edge corresponding to the most complex dynamical phase. Moreover, in the second half of the time interval the red line increases significantly. This could be due to an accumulation of the error, which is reasonable since the system is studied over a very large time interval. Let us try and consider more accurate integration methods, such as BDF3 of order $p=3$ to propagate the particle ensemble, and BFD4 to get an estimate of the error (Figure 2).

\begin{figure}
	\begin{center}
		\includegraphics[width=70mm]{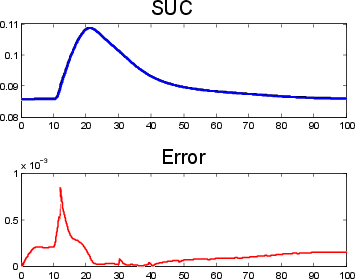}\\
		\includegraphics[width=70mm]{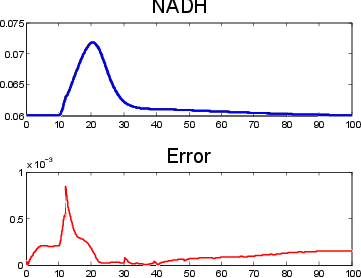}\\
		\includegraphics[width=70mm]{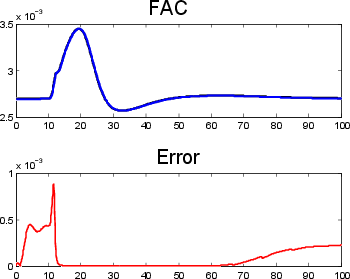}
	\end{center}
	\caption{Computed concentrations of succinate, fatty acyl-CoA and NAD reduced compared with their respective errors (BDF3-BDF4).}
\end{figure}

Although the error significantly decreases, it presents an increasing tendency towards the end of the time interval. Enlarging the observation time interval, it could be possible to observe a progressive separation between the blue and the black line, which is not reasonable comparing to the smoothness of the dynamic. It is worth emphasizing that the real solution of the ODEs system governing the problem is not known. The blue solid line is obtained by integrating the system with \verb|ode15s|, which is a variable order solver usually used to deal with stiff problems, such as in this case. Therefore, the accumulation of the error could be due to the usage of three numerical solvers.
\\To investigate the nature of the error, we apply the LMM PF to a test problem whose analytic solution is known. This allows us to avoid a possible source of error. We start by considering a very smooth problem:
\begin{equation}
x'(t)=cos^{2}(t),\quad x(0)=1,
\end{equation}
where 
\begin{equation*}
x(t)=arctg(t).
\end{equation*}
We need to set the number \verb|Nsample| of particles, the variance \verb|V| of the initial set of particles and the discretization step \verb|dt|:
\begin{verbatim}
Nsample=150;
V=0.1;
dt=0.1;
\end{verbatim}
We choose the Adams-Bashforth method of order $p=1$ to propagate the particle ensemble, and the Adams-Bashforth method of order $p=2$ to get an estimate of the error.

\begin{figure}
	\begin{center}
		\includegraphics[width=80mm]{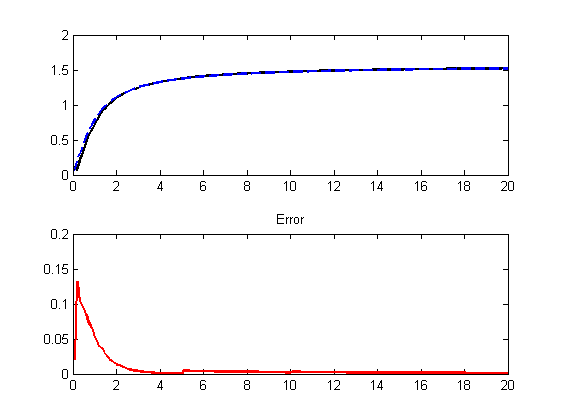}
	\end{center}
	\caption{Computed result of (11) compared with the respective error (AB1-AB2).}
\end{figure}

The error curve presents a reasonable edge corresponding to the increasing dynamical phase. Then, it does not increase anymore. The LMM PF algorithm does not worsen the stability of the numerical method chosen.
\\Let us consider a less smooth problem, such as
\begin{equation}
x'(t)=-2(t-1)x,\quad x(0)=1, \quad t\in [0,5],
\end{equation}
where 
\begin{equation*}
x(t)=e^{-t(t-2)}.
\end{equation*}
We are going to compute the solution for \verb|V=0.1,0.01,0.001,0.0001| when different integration methods are considered. The numerical results are compared in terms of the absolute error and of the sample variance at each time instant. In the following experiments the discretization step and the number of particles are fixed as before:
\begin{verbatim}
Nsample=150;
dt=0.1;
\end{verbatim}
We considered explicit and implicit linear multistep methods, such as Adams-Bashforth methods and Adams-Moulton methods. In particular, the adopted pair are AB1-AB2, AB3-AB4, AM1-AM2, AM3-AM4.

In Figure 4 and Figure 5, numerical results obtained with AB12 and AB34 respectively for \verb|V=0.1| and \verb|V=0.0001| are compared. 

\begin{figure}
	\begin{center}
		\includegraphics[width=70mm]{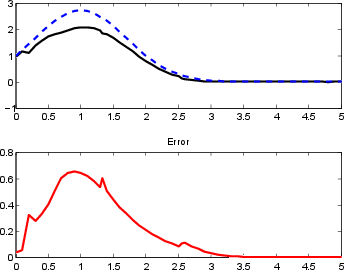}
		\includegraphics[width=70mm]{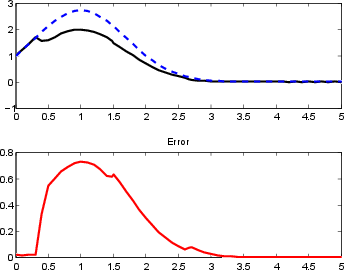}
	\end{center}
	\caption{Numerical results computed with AB12 (left) and AB34 (right) for $V=0.1$.}
\end{figure}

\begin{figure}
	\begin{center}
		\includegraphics[width=70mm]{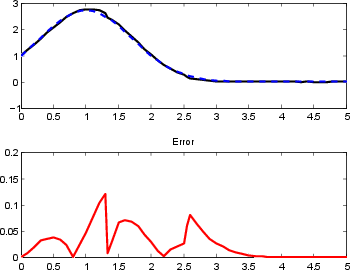}
		\includegraphics[width=70mm]{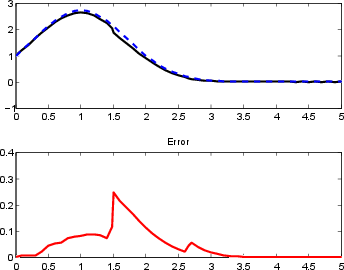}
	\end{center}
	\caption{Numerical results computed with AB12 (left) and AB34 (right) for $V=0.0001$.}
\end{figure}

Observe that the choice of a higher order method does not lead to a greater accuracy in the final solution. This is due to a reduction of the stability regions. An improvement could be observed decreasing $h$. The same issue occurs when considering the implicit pairs AM1-AM2 and AM3-AM4 (Figure 6, Figure 7).

\begin{figure}
	\begin{center}
		\includegraphics[width=70mm]{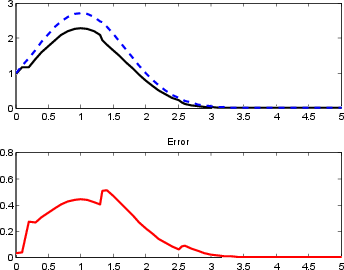}
		\includegraphics[width=70mm]{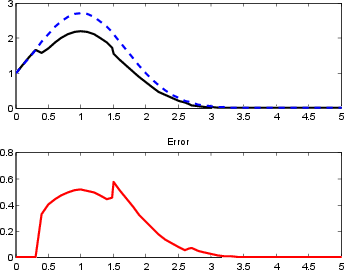}
	\end{center}
	\caption{Numerical results computed with AM12 (left) and AM34 (right) for $V=0.1$.}
\end{figure}

\begin{figure}
	\begin{center}
		\includegraphics[width=70mm]{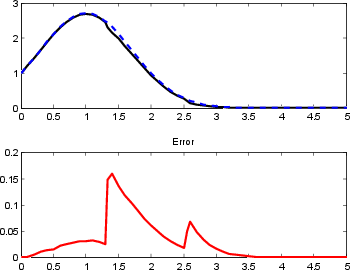}
		\includegraphics[width=70mm]{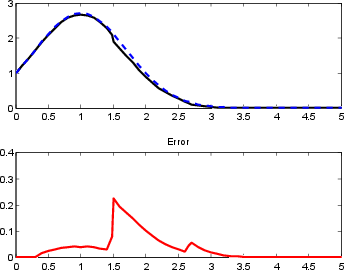}
	\end{center}
	\caption{Numerical results computed with AM12 (left) and AM34 (right) for $V=0.0001$.}
\end{figure}
Table 1-4 emphasize the worsening in the accuracy, by showing the norm of the absolute error vector and of the variance vector for the above mentioned linear multistep method when different values of the initial variance $V$ are chosen. Furthermore, as we would expect, the choice of implicit methods rather than the explicit ones ensures an higher accuracy.

\begin{table}
	\begin{tabular}{lp{0.3\textwidth}ll}
		Method & Absolute error ($\lVert\cdot\lVert_{\infty}$) & Sample variance ($\lVert\cdot\lVert_{2}$) \\
		$AB1-AB2$ & $0,4292$ & $0,1015$ &\\
		$AB3-AB4$ & $0,6795$ & $0,1011$  &\\ 
		$AM1-AM2$ & $0,2265$ & $0,1011$ &\\
		$AM3-AM4$ & $0,5728$ & $0,1005$ &
		
	\end{tabular}
	\caption{Numerical results for V=0.1}
\end{table}
As the initial variance V decreases, the distance of the generic particle of the sample from the sample mean decreases too, i.e. the LMM PF is more precise. 
\begin{table}
	\begin{tabular}{lp{0.3\textwidth}ll}
		Method (V=0.01) & Absolute error ($\lVert\cdot\lVert_{\infty}$) & Sample variance ($\lVert\cdot\lVert_{2}$) \\
		$AB1-AB2$ & $0,3724$ & $0,0218$ &\\
		$AB3-AB4$ & $0,6142$ & $0,0127$  &\\ 
		$AM1-AM2$ & $0,1545$ & $0,0181$ &\\
		$AM3-AM4$ & $0,5030$ & $0,0172$ & 
	\end{tabular}
	\caption{Numerical results for V=0.01}
\end{table}

\begin{table}
	\begin{tabular}{lp{0.3\textwidth}ll}
		Method (V=0.001) & Absolute error ($\lVert\cdot\lVert_{\infty}$) & Sample variance ($\lVert\cdot\lVert_{2}$) \\
		$AB1-AB2$ & $0,1216$ & $0,0107$ &\\
		$AB3-AB4$ & $0,3645$ & $0,0067$  &\\ 
		$AM1-AM2$ & $0,1132$ & $0,0061$ &\\
		$AM3-AM4$ & $0,2873$ & $0,0044$ &
	\end{tabular}
		\caption{Numerical results for V=0.001}
\end{table}

\begin{table}
	\begin{tabular}{lp{0.3\textwidth}ll}
		Method (V=0.0001)& Absolute error ($\lVert\cdot\lVert_{\infty}$) & Sample variance ($\lVert\cdot\lVert_{2}$) \\
		$AB1-AB2$ & $0,1140$ & $0,0059$ &\\
		$AB3-AB4$ & $0,2484$ & $0,0031$  &\\ 
		$AM1-AM2$ & $0,0681$ & $0,0014$ &\\
		$AM3-AM4$ & $0,2282$ & $0,0017$ &
	\end{tabular}
		\caption{Numerical results for V=0.0001}
\end{table}

Let us now test the Runge-Kutta integration class method on (12). It is worth remarking that, as the order of the RK method increases, the stability region enlarges. Hence, we can consider RK methods of order $p\geq3$ without worrying about the stability properties of the final solution and we can consider a larger integration time interval, such as $[0,10]$. In the following, numerical results obtained with RK12 and RK45 for different values of $V$ are compared.
\begin{figure}
	\begin{center}
		\includegraphics[width=70mm]{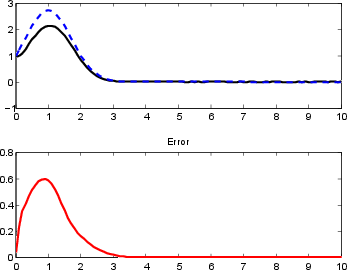}
		\includegraphics[width=70mm]{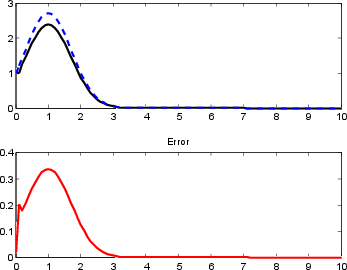}
	\end{center}
	\caption{Numerical results computed with RK12 (left) and RK45 (right) for $V=0.1$.}
\end{figure}

\begin{figure}
	\begin{center}
		\includegraphics[width=70mm]{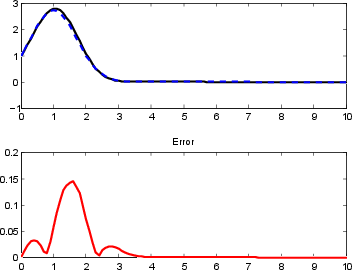}
		\includegraphics[width=70mm]{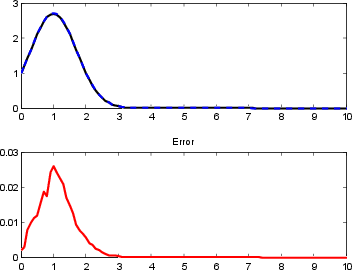}
	\end{center}
	\caption{Numerical results computed with RK12 (left) and RK45 (right) for $V=0.0001$.}
\end{figure}

Let us emphasize that a better accuracy has been obtained without changing the integration step \verb|dt=0.1|, that means that the integration time does not significantly increases. Moreover, fixing an admissible upper bound for the error, we can choose a not very accurate solver and a low initial variance, or a more accurate solver and a higher initial variance, i.e. we do not need to sample as better as possible at the initial time $t_{0}$.\\

\begin{table}
	\begin{tabular}{lp{0.3\textwidth}ll}
		Method (V=0.1)& Absolute error ($\lVert\cdot\lVert_{\infty}$) & Sample variance ($\lVert\cdot\lVert_{2}$) \\
		$RK1-RK2$ & $0,3337$ & $0,1024$ &\\
		$RK4-RK5$ & $0,2596$ & $0,1014$  &
	\end{tabular}
	\caption{Numerical results for V=0.1}
\end{table}

\begin{table}
	\begin{tabular}{lp{0.3\textwidth}ll}
		Method (V=0.01)& Absolute error ($\lVert\cdot\lVert_{\infty}$) & Sample variance ($\lVert\cdot\lVert_{2}$) \\
			$RK1-RK2$ & $0,3303$ & $0,0180$ &\\
			$RK4-RK5$ & $0,2862$ & $0,0152$  &
	\end{tabular}
	\caption{Numerical results for V=0.01}
\end{table}

\begin{table}
	\begin{tabular}{lp{0.3\textwidth}ll}
		Method (V=0.001)& Absolute error ($\lVert\cdot\lVert_{\infty}$) & Sample variance ($\lVert\cdot\lVert_{2}$) \\
			$RK1-RK2$ & $0,1574$ & $0,0088$ &\\
			$RK4-RK5$ & $0,1488$ & $0,0077$  &
	\end{tabular}
	\caption{Numerical results for V=0.001}
\end{table}

\begin{table}
	\begin{tabular}{lp{0.3\textwidth}ll}
		Method (V=0.0001)& Absolute error ($\lVert\cdot\lVert_{\infty}$) & Sample variance ($\lVert\cdot\lVert_{2}$) \\
			$RK1-RK2$ & $0,1418$ & $0,0035$ &\\
			$RK4-RK5$ & $0,0267$ & $0,0023$  &
	\end{tabular}
	\caption{Numerical results for V=0.0001}
\end{table}

\end{document}